\documentclass [12pt] {amsart}

\usepackage [T1] {fontenc}

\usepackage {amssymb}
\usepackage {amsmath,bm}
\usepackage {amsthm}
\usepackage{tikz}
\usepackage{geometry}

\usepackage{url}

 \usepackage[titletoc,toc,title]{appendix}
 
 \usepackage{tikz-cd}

\usetikzlibrary{arrows,decorations.pathmorphing,backgrounds,positioning,fit,petri}

\usetikzlibrary{matrix}

\usepackage{enumitem}
\usepackage{hyperref}
\hypersetup{
    colorlinks=true,
    linkcolor=blue,
    filecolor=magenta,      
    urlcolor=cyan,
    citecolor=magenta
}

\DeclareMathOperator {\td} {td}

\DeclareMathOperator {\pr} {pr}

\DeclareMathOperator {\SL} {SL}

\DeclareMathOperator {\GL} {GL}

\DeclareMathOperator {\an} {an}
\DeclareMathOperator {\SAtyp} {SAtyp}

\DeclareMathOperator {\Exp} {Exp}

\DeclareMathOperator {\alg} {alg}

\DeclareMathOperator {\rk} {rk}

\DeclareMathOperator {\Zcl} {Zcl}

\DeclareMathOperator {\Atyp} {Atyp}

\DeclareMathOperator {\C} {\mathbb{C}}
\DeclareMathOperator {\R} {\mathbb{R}}
\DeclareMathOperator {\h} {\mathbb{H}}
\DeclareMathOperator {\Z} {\mathbb{Z}}
\DeclareMathOperator {\Q} {\mathbb{Q}}

\DeclareMathOperator {\ppr} {Pr}

\DeclareMathOperator {\rar} {\rightarrow}
\DeclareMathOperator {\seq} {\subseteq}

\DeclareMathAlphabet\urwscr{U}{urwchancal}{m}{n}
\DeclareMathAlphabet\rsfscr{U}{rsfso}{m}{n}
\DeclareMathAlphabet\euscr{U}{eus}{m}{n}
\DeclareFontEncoding{LS2}{}{}
\DeclareFontSubstitution{LS2}{stix}{m}{n}
\DeclareMathAlphabet\stixcal{LS2}{stixcal}{m} {n}

\newcommand{\B}{\mathbb{B}}
\newcommand{\N}{\mathbb{N}}
\newcommand{\Bexp}{\mathbb{B}_{\exp}}
\newcommand{\Cexp}{\mathbb{C}_{\exp}}

\newcommand{\glpq}{\GL_2^+(\Q)}
\newcommand{\Qalg}{\Q^{\alg}}
\newcommand{\zbar}{\bar{z}}
\newcommand{\xbar}{\bar{x}}
\newcommand{\ybar}{\bar{y}}
\newcommand{\D}{\mathrm{D}}

\newcommand{\FAtyp}{\mathrm{FAtyp}}
\newcommand{\SFAtyp}{\mathrm{SFAtyp}}
\newcommand{\im}{\mathrm{Im}}

\newcommand{\dangle}[1]{\langle \langle #1 \rangle \rangle}
\newcommand{\sangle}[1]{\langle  #1 \rangle }

\theoremstyle {definition}
\newtheorem {definition}{Definition} [section]
\newtheorem {remark} [definition] {Remark}

\newtheorem{example} [definition] {Example}

\newtheorem* {notation} {Notation}

\theoremstyle {plain}

\newtheorem {theorem} [definition] {Theorem}
\newtheorem {proposition} [definition] {Proposition}

\newtheorem {conjecture} [definition] {Conjecture}

\theoremstyle {remark}

\calclayout

\begin {document}

\title[Existential Closedness and Zilber-Pink]{The Existential Closedness and Zilber-Pink Conjectures}

\author{Vahagn Aslanyan}
\address{School of Mathematics, University of Leeds, Leeds LS2 9JT, UK}
\curraddr{Department of Mathematics, University of Manchester, Manchester M13 9PL, UK}

\email{vahagn.aslanyan@manchester.ac.uk}

\thanks{Supported by Leverhulme Trust Early Career Fellowship ECF-2022-082 (at the University of Leeds) and partially supported by EPSRC Open Fellowship EP/X009823/1 and DKO Fellowship (at the University of Manchester). }

\date{14 March 2024}



\vspace*{-0.7cm}

\begin{abstract}
In this paper we survey the history of, and recent developments on, two major conjectures originating in Zilber's model-theoretic work on complex exponentiation -- Existential Closedness and Zilber-Pink. The main focus is on the modular versions of these conjectures and specifically on novel variants incorporating the derivatives of modular functions. The functional analogues of all the conjectures that we consider are theorems which are presented too. The paper also contains some new results and conjectures.
\end{abstract}

\dedicatory{To Boris Zilber on the occasion of his 75th birthday}

\maketitle
\tableofcontents

\vspace*{1cm}

\section{Introduction}

In the early 2000s Boris Zilber produced an influential body of work around the model theory of the complex exponential field $\C_{\exp}:=(\C; +, \cdot, \exp)$ where $\exp: z\to e^z$ is the exponential function \cite{Zilb-pseudoexp,Zilb-exp-sum-published,ZilbExpSum}. He showed that \textit{Schanuel's conjecture} (SC for short) on the transcendence properties of $\exp$  (see \S \ref{sec:Schanuel}) plays a central role in the model-theoretic properties of $\C_{\exp}$. However, the conjecture is out of reach -- it implies the algebraic independence of $e$ and $\pi$ over the rationals, which is a long-standing unsolved problem. This makes it hard to understand the model theory of $\C_{\exp}$. So Zilber constructed algebraically closed fields of characteristic $0$ equipped with a unary function, which satisfies some of the basic properties of $\C_{\exp}$ and, most importantly, (the analogue of) Schanuel's conjecture. He then isolated and axiomatised the ``most'' existentially closed ones among these exponential fields by a Hrushovski style amalgamation-with-predimension construction. These are called \textit{pseudo-exponential fields}. While these models are not existentially closed in the first-order sense, they are existentially closed in certain ``tame'' extensions. The axiom guaranteeing this is known as \textit{Strong Existential Closedness} or \textit{Strong Exponential Closedness}, or SEC for short. 

Zilber showed that his axiomatisation of pseudo-exponential fields is uncountably categorical. In particular, there is a unique pseudo-exponential field of cardinality continuum, denoted by $\B_{\exp}$. Zilber conjectured that $\B_{\exp}$ is isomorphic to $\C_{\exp}$. This is equivalent to the combination of two conjectures -- Schanuel's conjecture and the Strong Exponential Closedness conjecture stating that SEC holds on $\C_{\exp}$. A variant of the SEC conjecture, known as \textit{Exponential Closedness} or \textit{Existential Closedness} (EC for short) is currently an active research field in model theory. It states roughly that systems of equations involving field operations and the complex exponential function have solutions unless they are ``overdetermined'' (i.e. the number of independent equations is larger than the number of variables). The notion of overdetermined systems is in fact related to Schanuel's conjecture: a system is overdetermined if its solution would be a counterexample to Schanuel's conjecture. As the name suggests, SEC is a strong version of EC guaranteeing that under certain conditions systems of exponential equations have generic solutions.

Zilber's work on the model theory of complex exponentiation also gave rise to a Diophantine conjecture -- the \textit{Conjecture on Intersections with Tori}, or CIT for short. It states roughly that intersections of algebraic varieties with torsion cosets of algebraic tori, whose dimensions are larger than expected, are governed by finitely many torsion cosets of algebraic tori. The statement makes sense in the more general setting of semi-abelian varieties which gives rise to the \textit{Conjecture on Intersections with Semi-abelian Varieties}, or CIS for short. Both CIT and CIS were proposed by Zilber in \cite{Zilb-exp-sum-published} and independently by Bombieri, Masser, and Zannier in \cite{Bom-Mas-Zan}. The Manin-Mumford and Mordell-Lang conjectures are special cases of CIS. Zilber used CIT to deduce a uniform version of Schanuel's conjecture from itself which then was used to establish some partial results towards Exponential Closedness (see \cite{Zilb-exp-sum-published}). 

SC, EC, and CIT are quite general in form; replacing the exponential function by other transcendental functions often allows one to formulate analogues of these conjectures in other settings. Most notably, such analogues have been extensively explored for modular functions and, in particular, the $j$-invariant. However, these analogues are being studied for other reasons too: the modular analogue of Schanuel's conjecture is a special case of the Grothendieck-Andr\'e generalised period conjecture (see \cite[1.3 Corollaire]{bertolin}, \cite[\textsection 23.4.4]{andre}, and \cite[\S 6.3]{Aslanyan-Eterovic-Kirby-closure-op}), EC in that setting is a natural problem in complex geometry and model theory, and the analogue of CIT is a special case of the Zilber-Pink conjecture for (mixed) Shimura varieties, henceforth referred to as ZP. The latter was proposed by Pink (independently from Zilber and Bombieri, Masser, Zannier) as a far-reaching conjecture unifying the Andr\'e-Oort, Andr\'e-Pink-Zannier, Manin-Mumford, and Mordell-Lang conjectures \cite{Pink,Pink-2}. 

 Furthermore, from a model-theoretic point of view, if both SC and EC hold then in a certain sense they give a ``complete'' list of properties (non first-order axioms) of the function under consideration.\footnote{This means, in particular, that if EC holds then SC is the strongest possible transcendence statement about the function under consideration.} This is formalised by Zilber's categoricity and quasiminimality theorem in the exponential setting. There is no such theorem in the modular setting and there cannot be one, for the upper half-plane (hence the set of the reals) is definable from the graph of $j$, but the philosophy of SC and EC together giving a full description of the algebraic and transcendental properties of $j$ still applies. It is likely that a formal categoricity/quasiminimality result can be established for some relations defined in terms of $j$ (which give proper reducts of the complex field with $j$); this is part of our current research programme.

In this paper we present the above-mentioned conjectures in the exponential and modular settings, mostly focusing on the latter. As pointed out above, the modular variants of these conjectures are in part motivated by their exponential counterparts. However, there are some inherent differences between the two settings resulting in quite different methods and approaches, although some methods work in both contexts. One such difference is that unlike exponential functions, which are defined on the whole complex plane, modular functions are defined only on the upper half-plane. These spaces are ``geometrically different'', which accounts for different approaches to EC and ZP in these two settings. 
 This also makes the model-theoretic treatment of modular functions significantly harder. For example, direct counterparts of many aspects of Zilber's work on exponentiation, e.g.\ categoricity and quasiminimality, fail gravely in the modular setting (as explained above).

Further, modular functions satisfy third order differential equations as opposed to first order differential equations for exponential functions. So we can consider SC, EC, and ZP for modular functions together with their first two derivatives (the third one being algebraic over these). This generalisation makes the problems more challenging, but it also gives a deeper insight into them  by providing a broader model-theoretic picture. Let us briefly discuss two more reasons to consider SC, EC, and ZP for modular functions together with derivatives. Often when dealing with variants of these conjectures, not least their differential versions, even when derivatives are not considered, the approaches and techniques require looking at the derivatives anyway (see, for instance, \cite{Aslanyan-Eterovic-Kirby-Diff-EC-j,Aslanyan-weakMZPD}). Also, modular forms of weight 2 are the derivatives of modular functions (weight 0), which means that studying these problems for modular forms of weight 2 (without derivatives) is the same as studying them for the first derivatives of modular functions. 

We state several versions of the conjectures in this new setting, some of which have appeared in the literature but some are new. We then explain the relationship between these various conjectures and present their functional variants all of which were proven in recent years, save for Ax's original theorem proven in 1971.

\subsection{Abbreviations}
In the paper we will consider several variants of three conjectures -- Schanuel's conjecture, the Existential Closedness conjecture, and the Zilber-Pink conjecture. We use abbreviations to refer to those conjectures, and for the convenience of the reader we list some of these abbreviations below.

\begin{itemize}[leftmargin = *]

    \item[] \textbf{Schanuel}

    \item SC -- Schanuel Conjecture

    \item MSC -- Modular Schanuel Conjecture

    \item MSCD -- Modular Schanuel Conjecture with Derivatives

    \item[] \textbf{Existential Closedness }

        \item EC -- Existential Closedness or Exponential 
    Closedness

    \item MEC -- Modular Existential Closedness

    \item MECD -- Modular Existential Closedness with Derivatives

    \item[] \textbf{Zilber-Pink}
    \item CIT -- Conjecture on Intersections with Tori

    \item ZP -- Zilber-Pink

    \item MZP -- Modular Zilber-Pink

    \item MZPD -- Modular Zilber-Pink with Derivatives

\end{itemize}

\subsection{Dedication}

This paper is dedicated to Boris Zilber on the occasion of his 75th birthday, and is motivated by his work.  Boris was my DPhil supervisor (jointly with Jonathan Pila) at the University of Oxford from 2013 to 2017. His guidance has been instrumental in shaping my mathematical thinking and research interests, and his continued support, both throughout my DPhil and after that, has been tremendously helpful in my mathematical career. The hours spent with Boris at the Mathematical Institute and at Merton are some of my fondest memories of Oxford. I would like to thank him for everything and wish him a happy 75th birthday.

\section{The exponential setting}

In this section we look briefly at Zilber's work on model theory of complex exponentiation and the conjectures it gave rise to.

\subsection{Schanuel's conjecture and Exponential Closedness} \label{sec:Schanuel}

We begin by formulating Schanuel's conjecture.

\begin{conjecture}[Schanuel's Conjecture, SC {\cite[p. 30]{Lang-tr}}]
    For any $\mathbb{Q}$-linearly independent complex numbers $z_1,\ldots,z_n$
\[\td_\mathbb{Q}\mathbb{Q}(z_1,\ldots,z_n, e^{z_1}, \ldots, e^{z_n}) \geq n\] where $\td$ stands for transcendence degree.
\end{conjecture}

This conjecture is believed to capture \emph{all} transcendence properties of the exponential function. This can and will shortly be explained in a more precise sense. For now let us mention that Schanuel's conjecture for $n=2$ already implies the algebraic independence of $e$ and $\pi$ by choosing $z_1=\pi i,~ z_2=1$, which is a long-standing open problem. Thus, even for $n=2$ the conjecture is out of reach of current methods. However, partial results towards this conjecture are known, including the Lindemann-Weierstrass theorem and the Gelfond-Shneider theorem.

In \cite{Zilb-pseudoexp} Zilber presented a novel model-theoretic approach to Schanuel's conjecture. He constructed algebraically closed fields of characteristic $0$ equipped with a unary function, known as \textit{pseudo-exponentiation}, satisfying certain basic properties of the complex exponential functions and some desirable properties, not least the analogue of Schanuel's conjecture. He axiomatised these structures in the language $L_{\omega_1, \omega}(Q)$, where $Q$ is a quantifier for ``there are uncountably many'', and showed that the resulting theory is categorical in uncountable cardinals. The unique model of cardinality $2^{\aleph_0}$ is called \textit{the pseudo-exponential field} or \textit{the Zilber field} and is usually denoted by $\B_{\exp}$. Zilber then conjectured that $\B_{\exp}$ is isomorphic to $\C_{\exp}$. This shows, in a sense, that Schanuel's conjecture must play a central role in the model theory of $\C_{\exp}$. 

Zilber verified that all of the axioms of pseudo-exponentiation hold in $\C_{\exp}$ save for Schanuel's conjecture and an axiom called \textit{Strong Exponential Closedness} (SEC). So Zilber's conjecture that $\Bexp \cong \Cexp$ is equivalent to the conjunction of Schanuel's conjecture and the Strong Exponential Closedness conjecture (stating that the axiom holds in $\Cexp$). 

Let us explain what (Strong) Exponential Closedness means. Schanuel's conjecture can be interpreted as a statement about non-solvability of certain systems of equations which we demonstrate on an example below.
\begin{example}[\cite{Aslanyan-Kirby-Mantova}]
    Assume $e$ and $\pi$ are algebraically independent over $\Q$. Then for any non-constant polynomial $p(X,Y)\in \Q[X,Y]$ the system $e^z = 1,~ p(z, e) = 0$ does not have solutions in $\C$. On the other hand, if $e$ and $\pi$ are algebraically dependent, then for some $p$ that system does have a complex solution.
\end{example}

Another reason for a system not to have a solution is when the system is incompatible with the identity $e^{x+y}=e^xe^y$.

\begin{example}
    The system $z_2=z_1+1,~ 3e^{z_1}=e^{z_2}$ does not have a solution, for the first equation implies $e^{z_2} = e\cdot e^{z_1}$ and $e \neq 3$. On the other hand, the system $z_2=z_1+1,~ e^{z_2} = z_1,~ e^{z_2} = e\cdot e^{z_1}$ does have solutions even though there are three equations in two variables. Of course, the three equations are not ``analytically'' independent -- the third one follows from the first one by taking exponentials of both sides -- but they are algebraically independent.
\end{example}

In general, systems incompatible with the functional equation of $\exp$ are not solvable. Moreover, SC implies that if a system is ``overdetermined'', e.g. $n$ variables with more than $n$ algebraically independent equations, then there is no solution, unless the system can somehow be reduced using the functional equation $e^{x+y}=e^xe^y$.  With this interpretation SC becomes more natural, and Exponential Closedness (EC) is its dual conjecture stating roughly that a system of exponential equations does have a solution in $\C$ unless having a solution contradicts Schanuel's conjecture. Let us give a precise statement in geometric terms, observing first that  understanding the solvability of systems of exponential equations is equivalent to understanding when algebraic varieties contain exponential points (i.e. points of the form $(\zbar, \exp(\zbar))$). For instance, the equation $e^{e^z}+z-1=0$ has a solution if and only if the variety $V \seq \C^2 \times (\C^{\times})^2$ (with coordinates $(x_1,x_2,y_1,y_2)$ defined by the equations $x_2=y_1,~ y_2+x_1-1=0$ contains an exponential point.

\begin{conjecture}[Exponential Closedness, EC \cite{Zilb-pseudoexp,Bays-Kirby-exp}]

Let $V \seq \C^n \times (\C^{\times})^n$ be a free and rotund variety. Then $V$ contains a point of the form $(z_1,\ldots,z_n, e^{z_1}, \ldots, e^{z_n})$.
    
\end{conjecture}

 Freeness and rotundity are the conditions that make sure containing an exponential point does not contradict SC, as illustrated on the above examples. Now we define these notions precisely.

\begin{definition}
    An irreducible variety $V \subseteq \C^n \times (\C^{\times})^n$ is \emph{additively (resp. multiplicatively) free} if its projection to $\C^n$ (resp. $(\C^{\times})^n$)  is not contained in a translate of a $\Q$-linear subspace of $\C^n$ (resp. algebraic subgroup of $(\C^{\times})^n$). A variety is called \emph{free} if it is additively and multiplicatively free.
\end{definition}

We let $\xbar$ and $\ybar$ denote the coordinates on $\C^n$ and $(\C^{\times})^n$ respectively.  For a $k \times n$ matrix $M$ of integers we define $[M]:\C^n \times (\C^{\times})^n \rightarrow \C^k \times (\C^{\times})^k$ to be the map given by $[M]:(\bar{x},\bar{y}) \mapsto \left(M \xbar, \ybar^M\right)$ where
$$(M\xbar)_i = \sum_{j=1}^n m_{ij}x_j \mbox{ and } \left(\ybar^M\right)_i = \prod_{j=1}^n y_j^{m_{ij}}.$$

\begin{definition}\label{rotund}
An irreducible variety $V \subseteq \C^n \times (\C^{\times})^n$ is \emph{rotund} if for any $1 \leq k \leq n$ and any $k\times n$ matrix $M$ of integers $\dim [M](V) \geq \rk M$.
\end{definition}

Since $\exp$ maps $\Q$-linear equations to multiplicative ones, if the projections of $V$ satisfy either a linear or multiplicative equation and we want it to contain an exponential point, then these equations should match, otherwise they will not be compatible with $\exp$. Freeness takes care of this scenario by ensuring no such equations hold on the variety. Rotundity comes from SC; it states that $V$ and its various projections given by the maps $[M]$ have sufficiently large dimension so an exponential point in $V$ would not give a counterexample to SC.

Now we can formulate SEC which is a strong version of EC.

\begin{conjecture}[Strong Exponential Closedness, SEC \cite{Zilb-pseudoexp,Bays-Kirby-exp}]
    
Let $V \seq \C^n \times (\C^{\times})^n$ be a free and rotund variety. Then for every finitely generated field $K \seq \C$ over which $V$ is defined, there is a point $(z_1,\ldots,z_n, e^{z_1}, \ldots, e^{z_n})\in V$ which is generic in $V$ over $K$, that is, $\td_KK(\zbar, e^{\zbar}) = \dim V$.
\end{conjecture}

It is obvious that SEC implies EC. The converse is also true assuming SC and CIT hold (see \cite{eterovic-generic,Kirby-Zilber-exp}).

\begin{remark}
    The Rabinowitsch trick can be used to show that EC implies that a free and rotund variety contains a Zariski dense set of exponential points (see \cite[Theorem 4.11]{Kirby-semiab} and \cite[Proposition 4.34]{Aslanyan-adequate-predim}), but a priori such a set may not contain a generic point.
\end{remark}

\subsection{Conjecture on Intersections with Tori}

In \cite{Zilb-exp-sum-published} Zilber studied the solvability of exponential sums equations as a special case of the Exponential Closedness conjecture. In order to prove that certain systems of such equations are solvable, he needed a uniform version of Schanuel's conjecture. He then proposed a Diophantine conjecture, called the \emph{Conjecture on Intersections with Tori}, or \textit{CIT} for short, which acts as the difference between SC and Uniform SC. The conjecture states roughly that when we intersect an algebraic variety with algebraic tori then we do not expect to get too many intersections which are \textit{atypically large}. We will shortly give a precise formulation, but we need to introduce some notions first.

Let $V$ and $W$ be subvarieties of some variety $S$. A non-empty component $X$ of the intersection $V \cap W$ is \emph{atypical} in $S$ if $\dim X > \dim V + \dim W - \dim S$, and \emph{typical} if $\dim X = \dim V + \dim W - \dim S$. Note that if $S$ is smooth then a non-strict inequality always holds.

An \emph{algebraic torus} is an irreducible algebraic subgroup of a multiplicative group $(\C^{\times})^n$. Algebraic subgroups of $(\C^{\times})^n$ (not necessarily irreducible) are defined by multiplicative equations of the form $y_1^{m_1}\cdots y_n^{m_n}=1$ with $m_i \in \Z$. Any system of such equations (if consistent) defines an algebraic group. It splits as the union of an algebraic torus (the component containing the identity) and its finitely many translates by torsion points. Torsion cosets of algebraic tori are called \emph{special} varieties. For an algebraic variety $V \seq (\C^{\times})^n$ an \textit{atypical subvariety} of $V$ is an atypical component of an intersection of $V$ with a special variety $T$.  

Now we are ready to formulate the Conjecture on Intersections with Tori, which is the Zilber-Pink conjecture for algebraic tori. There are many equivalent forms of the conjecture; we consider four of them.

\begin{conjecture}[Conjecture on Intersections with Tori, CIT \cite{Zilb-exp-sum-published,Bom-Mas-Zan,Pila-ZP}]
Let $V \subseteq (\C^{\times})^n$ be an algebraic variety. 

\begin{enumerate}[leftmargin = *]
    \item[\rm{(1)}] There is a finite collection $\Sigma$ of proper special subvarieties of $(\C^{\times})^n$ such that every atypical subvariety of $V$ is contained in some $T\in \Sigma$.

    \item[\rm{(2)}] $V$ contains only finitely many maximal atypical subvarieties.

    \item[\rm{(3)}] Let $\Atyp(V)$ be the union of all atypical subvarieties of $V$. Then $\Atyp(V)$ is contained in a finite union of proper special subvarieties of $(\C^{\times})^n$.

    \item[\rm{(4)}] $\Atyp(V)$ is a Zariski closed subset of $V$.
\end{enumerate}

\end{conjecture}

\subsection{Functional/differential variants}

We have so far considered three conjectures for $(\C^{\times})^n$, namely, SC, EC, and CIT. As pointed out above, Schanuel's conjecture is out of reach, CIT is wide open and while EC is more tractable, it is also open. In spite of that, the functional analogues of all three conjectures are known.

Ax proved a functional analogue of Schanuel's conjecture in 1971. Below in a differential field $(F; +, \cdot, D_1, \ldots, D_m)$ we define a relation $\Exp(\xbar,\ybar)$ as the set of all $(\xbar, \ybar) \in F^n \times (F^{\times})^n$ for which $D_k y_i = y_i D_k x_i$ for all $k, i$.  Then $\Exp(F)$ is the set of all tuples $(\xbar,\ybar)\in F^n \times (F^{\times})^n$ with $F\models \Exp(\xbar, \ybar)$ (for all $n$).\footnote{ More generally, here and later, given a relation $R$ on a structure $M$, we write $R(M)$ for the set of all tuples from $M$ satisfying the relation $R$.}

\begin{theorem}[Ax-Schanuel {\cite[Thm. 3]{Ax}}]\label{thm:Ax-Schanuel-exp}
Let $(F;+,\cdot,D_1,\ldots,D_m)$ be a differential field with field of constants $C = \bigcap_{k=1}^m \ker D_k$. Let also $(x_i, y_i) \in F^2,~ i=1,\ldots,n,$ be such that $ (\xbar, \ybar)\in \Exp(F)$. Assume $x_1,\ldots,x_n$ are $\Q$-linearly independent mod $C$, that is, they are $\Q$-linearly independent in the quotient vector space $F/C$.
Then
$\td_CC\left(\bar{x},\bar{y}\right) \geq n+\rk (D_kx_i)_{i,k}.$
\end{theorem}

Ax's proof of this theorem is differential algebraic.  There is an equivalent complex analytic formulation of Ax-Schanuel (the equivalence follows from Seidenberg's embedding theorem). In \cite{Tsim} Tsimerman gave a new proof of that complex analytic statement based on o-minimality.

The differential version of EC for fields with several commuting derivations was established recently by Aslanyan, Eterovi\'c, and Kirby. 

\begin{theorem}[Differential EC {\cite[Thm. 4.3]{Aslanyan-Eterovic-Kirby-Diff-EC-j}}]
Let $(F; +, \cdot, D_1, \ldots, D_m)$ be a differential field with $m$ commuting derivations, and let $V\subseteq F^{2n}$ be a rotund variety. Then there exists a differential field extension $K$ of $F$ such that $V(K) \cap \Exp(K) \neq \emptyset$. In particular, when $F$ is differentially closed, $V(F) \cap \Exp(F) \neq \emptyset$.
\end{theorem}

The proof of this theorem uses some important differential algebraic ideas from \cite{Kirby-semiab} where the case of ordinary differential fields was treated. Kirby's approach (which in fact contains some inaccuracies and is not complete) is based on Ax's proof of the Ax-Schanuel theorem, while the argument given in \cite{Aslanyan-Eterovic-Kirby-Diff-EC-j} uses the statement of Ax-Schanuel as a black box and works quite generally.

\begin{example}\label{rmk:diff-ec-not-free}
    In the above theorem the variety $V$ need no be free. However, freeness is a necessary condition in EC. For example, the variety $V\seq \C^2 \times (\C^{\times})^2$ defined by the equations $x_2=x_1,~ y_2=y_1+1$, which is rotund but not free, cannot intersect the graph of any function. But it does intersect $\Exp(K)$ for any differential field $K$ -- indeed any constant point in $V$ is actually in $\Exp(K)$.
\end{example}

Finally, the following functional analogue of CIT was established independently by Zilber \cite{Zilb-exp-sum-published} and by Bombieri-Masser-Zannier \cite{Bom-Mas-Zan}. Both proofs rely on the Ax-Schanuel theorem. Kirby adapted Zilber's argument and gave a new proof in \cite{Kirby-semiab} using the uniform version of Ax-Schanuel, which follows from Ax-Schanuel by an application of the compactness theorem of first order logic (see \cite[Theorem 4.3]{Kirby-semiab}).

\begin{theorem}[Functional CIT \cite{Zilb-exp-sum-published,Bom-Mas-Zan,Kirby-semiab}]\label{thm:weak-CIT}
For every subvariety $V\subseteq (\C^{\times})^n$ there is a finite collection $\Sigma$ of proper subtori of $(\C^{\times})^n$ such that every atypical component of an intersection of $V$ with a coset of a torus is contained in a coset of some torus $T \in \Sigma$.
\end{theorem}

Theorem~\ref{thm:weak-CIT} is indeed a functional version of CIT as it talks about \textit{weakly special varieties} (arbitrary cosets of tori) and positive-dimensional atypical intersections.  In other words, it can be thought of as CIT over function fields where we work modulo the constants (in this case, the field of complex numbers). It does not say anything about special points or special coordinates in atypical intersections, so it is often called the geometric component of CIT (i.e. CIT without its arithmetic component). Since its statement is algebraic (rather than differential algebraic), it is also often called Weak CIT although, strictly speaking, it is not a weak version of CIT. 

In addition to the above-mentioned theorems, some other partial results have also been obtained towards EC and CIT in recent years. For EC see \cite{Aslanyan-Kirby-Mantova,Gallinaro-raising-to-powers,Gallinaro-exp-sums-trop,brown-masser,daquino--fornasiero-terzo}. It would be impractical to try to give a comprehensive list of references for CIT and its generalisations to semi-abelian varieties, so we refer the reader to Pila's recent book \cite{Pila-ZP} and references therein.

\section{The modular setting}

We let $\h\seq \C$ denote the complex upper half-plane and $j: \h \to Y(1)$ denote the modular $j$-function. We will identify the modular curve $Y(1)$ with the complex affine line $\C$.

Recall that the $j$-function is invariant under the linear fractional action of $\SL_2(\Z)$ on $\h$ and behaves nicely under the action of $\GL_2^+(\Q)$ (where $+$ denotes positive determinant). More precisely, there is a collection of \textit{modular polynomials} $\Phi_N(Y_1,Y_2)\in \Z[Y_1,Y_2],~N\in \N$, such that
\[ \text{for all } z_1, z_2 \in \h \left( \exists g \in \GL_2^+(\Q) \text{ with } z_2=gz_1 \text{ iff } \exists N\in \N \text{ such that } \Phi_N(j(z_1),j(z_2)) = 0 \right).\]

These correspondences are often referred to as the ``functional equations'' of the $j$-function. They are analogous to the functional equation $e^{x+y} = e^xe^y$ of the exponential function. This analogy allows one to state the modular counterparts of the exponential conjectures mentioned in the previous section, and that is what we do in this section. We focus on the $j$-function as other modular functions can be treated similarly, and often results about other modular functions can be deduced from those about $j$ since $j$ is a uniformiser for the modular group: it generates the field of all modular functions.

Now let us introduce some notation that will be used throughout the rest of the paper.

\begin{notation}
Let $n$ be a positive integer, $k \leq n$ and $1\leq i_1 < \ldots < i_k \leq n$. 
\begin{itemize}[leftmargin = *] 

\item Subsets of $\C^{2n}$ (e.g. $\h^n \times \C^n$) will be interpreted as subsets of $\C^n \times \C^n$, and we will denote the coordinates on this space by $(\xbar, \ybar)$.

\item $\ppr_{\xbar} : \C^{2n} \to \C^n$ is the projection to the first $n$ coordinates, and $\ppr_{\ybar} : \C^{2n} \to \C^n$ is the projection to the second $n$ coordinates.

    \item $\pr_{\bar{i}}:\C^{n} \rightarrow \C^{k}$ is the map
$\pr_{\bar{i}}:(x_1,\ldots,x_n)\mapsto (x_{i_1},\ldots,x_{i_k}).$

\item $\ppr_{\bar{i}}:\C^{2n}\rightarrow \C^{2k}$ denotes the map
$\ppr_{\bar{i}}:(\bar{x},\bar{y})\mapsto (\pr_{\bar{i}}\bar{x},\pr_{\bar{i}}\bar{y}).$

\item By abuse of notation we will let $j: \h^n \to \C^n$ denote all Cartesian powers of itself and $\Gamma_j\seq \h^n \times \C^n$ denote its graph.
\end{itemize}
\end{notation}

\subsection{Modular Schanuel conjecture and Modular Existential Closedness}

We begin by stating the analogue of Schanuel's conjecture for the $j$-function. It is a special case of the Grothendieck-Andr\'e generalised period conjecture (see \cite[1.3 Corollaire]{bertolin}, \cite[\textsection 23.4.4]{andre}, and \cite[\S 6.3]{Aslanyan-Eterovic-Kirby-closure-op}).

\begin{conjecture}[Modular Schanuel Conjecture, MSC]\label{MSC}
Let $z_1,\ldots,z_n \in \h$ be non-quadratic numbers with distinct $\GL_2^+(\Q)$-orbits. Then $$\td_{\Q}\Q(z_1,\ldots,z_n,j(z_1),\ldots,j(z_n))\geq n.$$
\end{conjecture}

Schneider's theorem, stating that if both $z$ and $j(z)$ are algebraic over $\Q$ then $z$ must be a quadratic irrational number, is a special case of this conjecture.

As in the exponential setting, this conjecture can be interpreted as a statement about non-solvability of certain systems of equations involving the $j$-function.

\begin{example}
Let $a, b \in \Q^{\alg}$ be algebraic \textit{non-special numbers}, that is,  their preimages under $j$ are not quadratic irrationals. By Schneider's theorem, these preimages cannot  be algebraic. Consider the system
\[
\begin{cases}
j(z_1) = a,\\
j(z_2) = b,\\
z_1^2 + z_2^2 +1 =0.
\end{cases}
\]
If this system has a solution, then $\td_{\Q}\Q(z_1,z_2,j(z_1),j(z_2)) = 1$. Hence by MSC, either $z_1$ or $z_2$ must be a quadratic number or they must be in the same $\GL_2^+(\Q)$-orbit. By our choice of $a$ and $b$, the numbers $z_1$ and $z_2$ are transcendental over $\Q$, hence non-quadratic. If they satisfy a relation $z_2=gz_1$ for some $g\in \GL_2^+(\Q)$ then, together with the equation $z_1^2+z_2^2+1= 0$, we can conclude that $z_1,z_2\in \Q^{\alg}$, which is a contradiction. So MSC implies that the above system has no complex solutions. Note that it is overdetermined in the sense that we have 3 equations but only 2 variables.
\end{example}

Thus, we can propose a dual conjecture stating roughly that such a system always has a solution unless it contradicts MSC. We begin by recalling a few definitions from \cite{Aslanyan-adequate-predim}.

\begin{definition}\label{defin: broad-free}
    Let $V \subseteq \h^n \times \C^{n}$ be an algebraic variety.
\begin{itemize}[leftmargin = *]

    \item $V$ is \emph{$\Gamma_j$-broad} if for any $1\leq k_1 < \ldots < k_l \leq n$ we have $\dim \Pr_{\bar{k}} V \geq l$. 
    
    \item $V$ is \emph{modularly free} if no equation of the form $\Phi_N(y_i,y_k)=0$, or of the form $y_i =c$ with $c\in \C$ a constant, holds on $V$.
    
    \item  $V$ is $\glpq$-\emph{free} if no equation of the form $x_i = gx_k$ with $g \in \GL_2^+(\Q)$, {or of the form $x_k = c$ with $c\in \h$ a constant,} holds on $V$.

    \item  $V$ is $\Gamma_j$-\emph{free} if it is $\glpq$-free and modularly free.

    \item $V$ is $\Gamma_j$-\emph{froad}\footnote{To be pronounced like ``fraud''.} if it is $\Gamma_j$-free and $\Gamma_j$-broad.
\end{itemize}

\end{definition}

Now we are ready to state the Existential Closedness conjecture.

\begin{conjecture}[Modular Existential Closedness, MEC, {\cite[Cnj. 1.2]{Aslanyan-Kirby-blurred-j}}]\label{j-ec}
Let $V \subseteq \h^n \times \C^n$ be an irreducible $\Gamma_j$-froad variety defined over $\C$. Then $V\cap \Gamma_j \neq \emptyset$.
\end{conjecture}

As in the exponential setting, we can consider a strong version of MEC -- referred to as SMEC -- stating that $\Gamma_j$-froad varieties contain generic points from $\Gamma_j$. Eterovi\'c \cite{eterovic-generic} proved that MSC, MZP (see below), MEC imply SMEC. 

\subsection{Modular Zilber-Pink}

In \cite{Pink,Pink-2} Pink proposed a far-reaching conjecture in the setting of mixed Shimura varieties generalising the Manin-Mumford, Mordell-Lang, and Andr\'e-Oort conjectures. That conjecture also generalises Zilber's CIT conjecture (although Pink came up with it independently from Zilber and Bombieri-Masser-Zannier) and is now known as the \textit{Zilber-Pink conjecture}. Thus, CIT is the Zilber-Pink conjecture for algebraic tori. In this section we look at the Zilber-Pink conjecture in the modular setting, i.e. for $Y(1)^n$ (identified with $\C^n$ as usual).

\begin{definition}

\begin{itemize}[leftmargin = *]
    \item[]

    \item A $j$-\emph{special} variety in $\C^n$ is an irreducible component of a variety defined by some modular equations $\Phi_N(y_k, y_l) = 0$.

    \item Let $V\seq \C^n$ be a variety. A $j$-\emph{atypical} subvariety of $V$ is an atypical component of an intersection $V\cap T$ where $T$ is $j$-special.
\end{itemize}
 
\end{definition}

As for CIT, Modular Zilber-Pink has several equivalent formulations. Four of them are presented below.

\begin{conjecture}[Modular Zilber-Pink, MZP {\cite[Cnj. 19.2]{Pila-ZP}}]
Let $V \subseteq \C^n$ be an algebraic variety. Let also $\Atyp_j(V)$ be the union of all $j$-atypical subvarieties of $V$.  Then the following equivalent conditions hold.

\begin{enumerate}[leftmargin = *]
    \item[\rm{(1)}] There is a finite collection $\Sigma$ of proper $j$-special subvarieties of $\C^n$ such that every $j$-atypical subvariety of $V$ is contained in some $T\in \Sigma$.

    \item[\rm{(2)}] $V$ contains only finitely many maximal $j$-atypical subvarieties.

    \item[\rm{(3)}]  $\Atyp_j(V)$ is contained in a finite union of proper $j$-special subvarieties of $\C^n$.

    \item[\rm{(4)}] $\Atyp_j(V)$ is a Zariski closed subset of $V$.
\end{enumerate}

\end{conjecture}

As in the exponential setting, MZP and SC imply a uniform version of SC.

\subsection{Functional/differential variants}\label{subsec:j-diff-versions}

The $j$-function satisfies an order $3$ algebraic differential equation over $\mathbb{Q}$, and none of lower order (see \cite{mahler}). Namely, $\Psi(j,j',j'',j''')=0$ where 
$$\Psi(y_0,y_1,y_2,y_3)=\frac{y_3}{y_1}-\frac{3}{2}\left( \frac{y_2}{y_1} \right)^2 + \frac{y_0^2-1968y_0+2654208}{2y_0^2(y_0-1728)^2}\cdot y_1^2.$$

Thus
$$\Psi(y,y',y'',y''')=Sy+R(y)(y')^2,$$
where $S$ denotes the \emph{Schwarzian derivative} defined by $Sy = \frac{y'''}{y'} - \frac{3}{2} \left( \frac{y''}{y'} \right) ^2$ and $R(y)=\frac{y^2-1968y+2654208}{2y^2(y-1728)^2}$ is a rational function.

All functions $j(gz)$ with $g \in \SL_2(\mathbb{C})$ satisfy the differential equation $\Psi(y,y',y'',y''')=0$ and all solutions (not necessarily defined on $\mathbb{H}$) are of that form (see {\cite[Lemma 4.2]{Freitag-Scanlon}}). 

Note that for non-constant $y$, the relation $\Psi(y,y',y'',y''') = 0$ is equivalent to $y''' = \eta(y,y',y'')$ where $$\eta(y,y',y'') := \frac{3}{2}\cdot \frac{(y'')^2}{y'} - R(y) \cdot (y')^3$$ is a rational function over $\mathbb{Q}.$

\textbf{From now on, $y', y'', y'''$ will denote some variables/coordinates and not the derivatives of $y$. Derivations of abstract differential fields will not be denoted by $'$. When we deal with actual functions though, $'$ will denote the derivative, e.g. $j'$ is the derivative of $j$.}

\begin{definition}
Let $(F; +, \cdot, D_1,\ldots,D_m)$ be a differential field with constant field $C = \bigcap_{k=1}^m \ker D_k$. We define a binary relation $\D \Gamma_j(x, y)$  by
$$ \exists y', y'', y''' \left[ \Psi(y,y',y'',y''')=0 \wedge \bigwedge_{k=1}^m D_ky=y'D_kx \wedge D_ky'=y''D_kx \wedge D_ky''=y'''D_kx\right].$$
The relation $\D \Gamma_j^{\times}(x,y)$ is defined by the formula
$\D \Gamma_j(x,y)\wedge x \notin C \wedge y \notin C.$  By abuse of notation, we let $\D \Gamma_j$ and $\D \Gamma_j^{\times}$ also denote the Cartesian powers of these relations.
\end{definition}

 If $F$ is a field of meromorphic functions of variables $t_1,\ldots,t_m$ over some complex domain with derivations $\frac{d}{d t_k}$, then $\D \Gamma_j^{\times}(F)$ is interpreted as the set of all tuples $(x,y)\in F^2$ where $x=x(t_1,\ldots,t_m)$ is some meromorphic function and $y=j(gx)$ for some $g\in \GL_2(\C)$.

Pila and Tsimerman proved the following analogue of Ax's theorem for the $j$-function.

\begin{theorem}[Ax-Schanuel for $j$, {\cite[Thm. 1.3]{Pila-Tsim-Ax-j}}]\label{j-chapter-Ax-for-j}
Let $(F;+,\cdot,D_1,\ldots,D_m)$ be a differential field with commuting derivations and with field of constants $C$. Let also $(z_i, j_i) \in \D \Gamma_j^{\times}(F),~ i=1,\ldots,n.$ If the $j_i$'s are pairwise modularly independent  (i.e. no two of them satisfy an equation given by a modular polynomial) then $\td_CC\left(\bar{z},\bar{j}\right) \geq n+\rk (D_kz_i)_{i,k}.$
\end{theorem}

The proof of Pila and Tsimerman relies on o-minimality and, in particular, the definability of the restriction of the $j$-function to a fundamental domain in the o-minimal structure $\R_{\an,\exp}$. Recently, a differential-algebraic proof of Ax-Schanuel for all Fuchsian automorphic functions (including $j$) was given in \cite{sanz-cas-frei-nag-Ax-I}.

In \cite{Aslanyan-Eterovic-Kirby-Diff-EC-j}, Aslanyan, Eterovi\'c, and Kirby use the Ax-Schanuel theorem for the $j$-function to establish an existential closedness result for $\D \Gamma_j$. The proof is differential algebraic, and its advantage is that it treats Ax-Schanuel as a black box without looking into it, as opposed to the approach of \cite{Kirby-semiab} where the proof of Ax-Schanuel is appealed to. For that reason the proof works both for $\exp$ and $j$, and is expected to work in any reasonable situation where Ax-Schanuel is known.

\begin{theorem}[Functional MEC {\cite[Thm. 1.1]{Aslanyan-Eterovic-Kirby-Diff-EC-j}}]\label{thm-ec-j}
Let $F$ be a differential field, and $V\subseteq F^{2n}$ be a $\Gamma_j$-broad variety. Then there is a differential field extension $K\supseteq F$ such that $V(K)\cap \D \Gamma_j(K) \neq \emptyset$. In particular, if $F$ is differentially closed then $V(F)\cap \D \Gamma_j(F) \neq \emptyset$.
\end{theorem}

\begin{remark}
    In the above theorem the variety $V$ need not be free. However, freeness is a necessary condition in MEC; see Example~\ref{rmk:diff-ec-not-free}.  
    
      Also, when $V$ is defined over the constants $C$ and is \textit{strongly $\Gamma_j$-broad} (i.e. strict inequalities hold in Definition \ref{defin: broad-free} (first bullet point)), we have $V(K)\cap \D\Gamma_j^{\times}(K) \neq \emptyset$; see \cite[Theorem 1.3]{Aslanyan-Eterovic-Kirby-Diff-EC-j}.
\end{remark}

The Ax-Schanuel theorem can also be used to establish a functional variant of Modular Zilber-Pink, which was done by Pila and Tsimerman \cite[Theorem 7.1]{Pila-Tsim-Ax-j}. They used tools of o-minimality, while in \cite[Theorem 5.2]{Aslanyan-weakMZPD} Aslanyan gave a differential-algebraic proof based on Kirby's adaptation of Zilber's proof of weak CIT (see \cite[Theorem 4.6]{Kirby-semiab}).

\begin{definition}
Let $V \subseteq \C^n$ be an algebraic variety. A \emph{$j$-atypical subvariety} of $V$  is an irreducible component $W$ of some $V \cap T$, where $T$ is a $j$-special variety, such that $\dim W> \dim V + \dim T -n.$
A $j$-atypical subvariety $W$ of $V$ is said to be \emph{strongly $j$-atypical} if no coordinate is constant on $W$.
\end{definition}

\begin{theorem}[Functional MZP \cite{Pila-Tsim-Ax-j,Aslanyan-weakMZPD}]\label{thm:func-MZP}
Every algebraic variety $V \seq \C^n$ contains only finitely many maximal strongly $j$-atypical subvarieties.
\end{theorem}

Like the MZP conjecture, this theorem can also be stated in several equivalent forms, but we do not present them. See \cite{Aslanyan-weakMZPD} for details.

As in the exponential setting, recent years have seen significant progress towards MEC and MZP. For the state-of-the-art on MZP and its generalisations see \cite{Pila-ZP} and references therein. For MEC the reader is referred to \cite{eterovic-herrero,Aslanyan-Kirby-blurred-j,Gallinaro-raising-to-powers,eterovic-generic,eterovic-zhao}.

\section{Incorporating the derivatives of modular functions}\label{sec: incorporating derivatives}

In this section we look at the extensions of MSC, MEC, and MZP to the $j$-function together with its derivatives. Analogues of MSC and MZP in this setting were considered by Pila in some unpublished notes \cite{Pila-MZPD}, and we closely follow him in \S \ref{subsec:MSCD} and the beginning of \S \ref{subsec:MZPD}. MSC with Derivatives is in fact a special case of the Grothendieck-Andr\'e generalised period conjecture. MEC with Derivatives was first proposed in \cite{Aslanyan-Kirby-blurred-j} by Aslanyan and Kirby. In addition to that conjecture we also propose a second, more general MEC with Derivatives conjecture here. 

Recall that $j$ satisfies a third order differential equation, hence it suffices to consider only the first two derivatives. Adding higher derivatives would not change the problems. One normally works in jet spaces when dealing with $j$ together with its derivatives $j', j''$. However, as usual, instead of the jet space $J_2\h^n\times J_2 Y(1)^n$ we will work in $\h^n \times \C^{3n}$. We will use $(\xbar, \ybar, \ybar', \ybar'')$ to denote the coordinates on this space. We denote the vector-function $(j,j',j''): \h^n \to \C^{3n},~ \zbar \mapsto (j(\zbar), j'(\zbar), j''(\zbar))$ by $J$, and its graph by $\Gamma_J$.

Before proceeding we introduce further notation to be used in the rest of this section.

\begin{notation}
Let $n$ be a positive integer, $k \leq n$ and $1\leq i_1 < \ldots < i_k \leq n$. 
\begin{itemize}[leftmargin = *] 

\item $\Pi_{\bar{i}}:\C^{4n}\rightarrow \C^{4k}$ is defined by
$\Pi_{\bar{i}}:(\bar{x},\bar{y},\bar{y}',\bar{y}'')\mapsto (\pr_{\bar{i}}\bar{x},\pr_{\bar{i}}\bar{y},\pr_{\bar{i}}\bar{y}',\pr_{\bar{i}}\bar{y}'').$

\item $\pi_{\bar{i}}:\C^{3n}\rightarrow \C^{3k}$ denotes the map
$\pi_{\bar{i}}:(\bar{y},\bar{y}',\bar{y}'')\mapsto (\pr_{\bar{i}}\bar{y},\pr_{\bar{i}}\bar{y}',\pr_{\bar{i}}\bar{y}'').$

\item $\pi_{\bar{y}}:\C^{3n}\rar \C^n$, $\Pi_{\bar{y}}:\C^{4n}\rar \C^n$, and $\Pi_{\bar{x}}:\C^{4n}\to \C^n$ are the maps $\pi_{\bar{y}}:  (\bar{y},\bar{y}',\bar{y}'')\mapsto \bar{y}, \Pi_{\bar{y}}:  (\bar{x},\bar{y},\bar{y}',\bar{y}'')\mapsto \bar{y}, \mbox{ and } \Pi_{\bar{x}}:  (\bar{x},\bar{y},\bar{y}',\bar{y}'')\mapsto \bar{x}$ respectively.
\end{itemize}
\end{notation}

\subsection{Modular Schanuel Conjecture with Derivatives}\label{subsec:MSCD}

\begin{conjecture}[Modular Schanuel Conjecture with Derivatives; MSCD]\label{conj:MSCD-1}
Let $z_1,\ldots,z_n \in \h$ be non-quadratic numbers with distinct $\glpq$-orbits. Then $$\td_{\Q}\Q(z_1,\ldots,z_n,J(z_1),\ldots,J(z_n))\geq 3n.$$
\end{conjecture}

This conjecture is a direct generalisation of MSC, but it does not reflect the transcendence properties of $J$ at special points. So, following Pila \cite{Pila-MZPD}, we formulate a more general conjecture.

\begin{definition}\label{defin-intro-J-special-0}
\begin{itemize}[leftmargin = *] 
    \item[] 
    
    \item An irreducible subvariety $U\subseteq \mathbb{H}^n$ (i.e. an intersection of $\mathbb{H}^n$ with some algebraic variety) is called $\GL_2^+(\Q)$-\emph{special} if it is defined by some equations of the form $z_i = g_{i,k}z_k,~ i\neq k$, with $g_{i,k} \in \GL_2^+(\mathbb{Q})$, or of the form $z_i = \tau_i$ where $\tau_i \in \mathbb{H}$ is a quadratic number.
    
    \item  For a $\GL_2^+(\Q)$-special variety $U$ we denote by $\sangle{U}$ the Zariski closure of the graph of the restriction $J|_U$ (i.e. the set $\{ (\zbar, J(\zbar)): \zbar \in U \}$) over $\Qalg$.

    \item The $\glpq$-\textit{special closure} of an irreducible variety $W\seq \h^n$ is the smallest $\glpq$-special variety containing $W$. It exists for the irreducible components of an intersection of $\glpq$-special varieties is $\glpq$-special.
\end{itemize}
\end{definition}

We now explain how $\sangle{U}$ can be defined algebraically. First let us ignore the case when $U$ has constant coordinates. Assume the first two coordinates of $U$ are related, i.e. $x_2=gx_1$ for some $g = \begin{pmatrix} a & b \\ c & d \end{pmatrix}\in \glpq$, and let $\Phi(j(z),j(gz))=0$ for some modular polynomial $\Phi$. Differentiating the last equality with respect to $z$ we get
\begin{equation}
    \frac{\partial \Phi}{\partial Y_1} (j(z),j(gz)) \cdot j'(z) +  \frac{\partial \Phi}{\partial Y_2} (j(z),j(gz)) \cdot j'(gz) \cdot \frac{ad-bc}{(cz+d)^2} = 0. \tag{$\star$} \label{eq:deriv-of-j}
\end{equation}
Thus, $\sangle{U}$ satisfies the following equation:
\[     \frac{\partial \Phi}{\partial Y_1} (y_1,y_2) \cdot y'_1 +  \frac{\partial \Phi}{\partial Y_2} (y_1,y_2) \cdot y'_2 \cdot \frac{ad-bc}{(cx_1+d)^2} = 0. \tag{$\dagger$}\label{eq:<U>} \]
Differentiating once more we will get another equation between $(x_1,x_2,y_1,y_2,y_1',y_2',y_1'',y_2'')$, and we will have four equations defining the projection of $\sangle{U}$ to the first two coordinates. 

 In general, we have a partition of $\{ 1, \ldots, n\}$ where two indices are in the same block of the partition if and only if the corresponding coordinates are related on $U$. If $i_1<\ldots<i_k$ form such a block, then $\Pi_{\bar{i}}\sangle{U}$ will be referred to as a \textit{block} of $\sangle{U}$. Then each block of $\sangle{U}$ is defined by equations of the form described above and has dimension $4$, and $\sangle{U}$ is the product of its blocks.

When $U$ has a constant coordinate, say $x_1$ (whose value must be a quadratic irrational), then we also get blocks of dimension $1$ or $0$ as follows. If $x_1=\tau \notin \SL_2(\Z)i\cup \SL_2(\Z)\rho$, where $\rho = -\frac{1}{2}+\frac{\sqrt{3}}{2}i$, then $j(\tau)\in \Qalg$ and $\td_{\Q}\Q(j'(\tau), j''(\tau)) = 1$ (see \cite{diaz}). If, in addition, some other coordinates, say $x_2, \ldots, x_k$ are $\glpq$-related to $x_1$ and thus take constant values $\tau_k$ (with $\tau_1:=\tau$) then $\td_{\Q}(\bar{\tau}, J(\bar{\tau})) = 1$. Thus, we get a block of dimension $1$. The equations defining such a block can be worked out as above.

On the other hand, a constant coordinate in $\SL_2(\Z)\rho$ would give rise to a block of dimension $0$, for the values of $j,j',j''$ are zeroes at these points. A constant coordinate in $\glpq \rho \setminus \SL_2(\Z)\rho$ or $\glpq i$ gives a block of dimension $1$.

Now we are ready to state the second (and more general) version of MSCD.

\begin{conjecture}[Modular Schanuel Conjecture with Derivatives and Special Points; MSCDS]\label{conj:MSCD-2}
Let $z_1,\ldots,z_n \in \h$ be arbitrary and let $U \seq \h^n$ be the $\GL_2^+(\Q)$-special closure of $(z_1,\ldots,z_n)$. Then $$\td_{\Q}\Q(z_1,\ldots,z_n,J(z_1),\ldots,J(z_n))\geq \dim \sangle{U} - \dim U.$$
\end{conjecture}

Both MSCD and MSCDS are special cases of the Grothendieck-Andr\'e generalised period conjecture; see \cite[\S 6.3]{Aslanyan-Eterovic-Kirby-closure-op}.

\subsection{Modular Existential Closedness with Derivatives}\label{subsec:MECD}

We now introduce the appropriate notions of broadness and freeness which will appear in Existential Closedness.

\begin{definition} \label{defin: J broad-free}

Let $V \subseteq \h^n \times \C^{3n}$ be an algebraic variety.

\begin{itemize}[leftmargin = *]

    \item $V$ is \emph{$\Gamma_J$-broad} if for any $1\leq i_1 < \ldots < i_k \leq n$ we have $\dim \Pi_{\bar{i}} (V) \geq 3k$.
    
    \item $V$ is \emph{modularly free} if no equation of the form $\Phi_N(y_i,y_k)=0$, or of the form $y_i =c$ with $c\in \C$ a constant, holds on $V$.
    
    \item  $V$ is $\glpq$-\emph{free} if no equation of the form $x_i = gx_k$ with $g \in \GL_2^+(\Q)$, {or of the form $x_k = c$ with $c\in \h$ a constant,} holds on $V$.

    \item  $V$ is $\Gamma_J$-\emph{free} if it is $\glpq$-free and modularly free.

    \item $V$ is $\Gamma_J$-\emph{froad} if it is $\Gamma_J$-free and $\Gamma_J$-broad.
\end{itemize}

\end{definition}

\begin{conjecture}[Modular Existential Closedness with Derivatives, MECD]\label{conj:ECD-1}
Let $V \subseteq \h^n \times \C^{3n}$ be a $\Gamma_J$-froad variety defined over $\C$. Then $V\cap \Gamma_J \neq \emptyset$.
\end{conjecture}

This is dual to Conjecture \ref{conj:MSCD-1}. It is possible to state a dual to Conjecture \ref{conj:MSCD-2} which would also imply that certain varieties contain $J$-special points. However, in that case only dimension conditions would not suffice to guarantee existence of $J$-points, e.g. an arbitrary variety of dimension $1$ may not contain such a point; it should be $J$-special in order to contain $J$-special points. So we give the following definitions.

\begin{definition}
   Let $V \seq \h^n \times \C^{3n}$ be an irreducible variety. Let also $U \seq \h^n$ be the $\GL_2^+(\Q)$-special closure of $\Pi_{\xbar} (V)$ and $T \seq \C^n$ be the $j$-special closure\footnote{The $j$-special closure of an irreducible set $W$ is defined as the smallest $j$-special set containing $W$.} of $\Pi_{\ybar}(V)$.

   \begin{itemize}[leftmargin = *]
       \item $V$ is said to be $\Gamma_J^*$-\textit{free} if $j(U)=T$ and $V\seq \sangle{U}$;

       \item $V$ is said to be $\Gamma_J^*$-\textit{broad} if for any $\bar{i}$ we have
   $\dim \Pi_{\bar{i}} (V) \geq \dim \sangle{\pr_{\bar{i}} U} - \dim \pr_{\bar{i}} U. $

    \item $V$ is said to be $\Gamma_J^*$-\textit{froad} if it is $\Gamma_J^*$-free and $\Gamma_J^*$-broad.
   \end{itemize}
\end{definition}

\begin{remark}
$\Gamma_J^*$-freeness means that the $\glpq$-relations and modular relations holding on $V$ match each other, i.e. are compatible with the functional equations of $J$ (that is, modular correspondences and the relations obtained by differentiating those). This condition holds vacuously for $\Gamma_J$-free varieties. For $\Gamma_J$-free varieties $\Gamma_J$-broadness and $\Gamma_J^*$-broadness are equivalent.
\end{remark}

\begin{conjecture}[Modular Existential Closedness with Derivatives and Special Points, MECDS]\label{conj:ECD-2}
    Let $V \subseteq \h^n \times \C^{3n}$ be an irreducible $\Gamma_J^*$-froad variety. Then $V\cap \Gamma_J \neq \emptyset$.
\end{conjecture}

\subsection{Modular Zilber-Pink with Derivatives}\label{subsec:MZPD}

\begin{definition}\label{defin-intro-J-special}
\begin{itemize}[leftmargin = *] 
    \item[]

    \item  For a $\GL_2^+(\Q)$-special variety $U\seq \h^n$ we denote by $\langle \langle U \rangle \rangle$ the Zariski closure of $J(U)$ over $\mathbb{Q}^{\alg}$.
    
    \item A $J$-\textit{special} subvariety of $\mathbb{C}^{3n}$ is a set of the form $S = \langle \langle U \rangle \rangle$ where $U$ is a $\glpq$-special subvariety of $\mathbb{H}^n$. 

    \item A $J$-special variety $S$ is said to be \textit{associated to }a $j$-special variety $T$ if there is a $\glpq$-special variety $U$ such that $S=\dangle{U}$ and $j(U)=T$.
\end{itemize}

\end{definition}

\begin{remark}\label{remark-image-analytic}
\begin{itemize}[leftmargin = *] 
    \item[]

    \item For a $\glpq$-special variety $U\subseteq \mathbb{H}^n$ the set $j(U)\subseteq \mathbb{C}^n$ is defined by modular equations and is irreducible (since $U$ is irreducible), therefore it is $j$-special. Similarly, $J(U)$ is an irreducible locally analytic set\footnote{Strictly speaking, $J(U)$ may not be complex analytic as it is the image of an analytic set under an analytic function, but it is locally analytic. It is irreducible in the sense that if $J(U)$ is contained in a countable union of analytic sets then it must be contained in one of them.} and hence so is its Zariski closure. Thus, $J$-special varieties are irreducible.

    \item $j$-special varieties are bi-algebraic for the $j$-function, that is, they are the images under $j$ of algebraic varieties (namely, $\glpq$-special varieties). That is in contrast to $J$-special varieties as these are not bi-algebraic for $J$. Nonetheless, $J$-special varieties still capture the algebraic properties of the function $J$.

    \item The equations defining a $J$-special variety can be worked out as in \S \ref{subsec:MSCD} since $\dangle{U}$ is a projection of $\sangle{U}$. In particular, a variety $\dangle{U}$ is the product of its \textit{blocks} each of which has dimension $0, 1, 3$ or $4$. Dimensions 0 and 1 correspond to constant coordinates. A block has dimension 3 if all the $\glpq$-matrices linking its $x$-coordinates are upper triangular, and dimension 4 otherwise. This is because equation \eqref{eq:<U>} gives an algebraic relation between $y_1,y_2,y_1',y_2'$ provided that $c=0$, i.e. the matrix linking $x_1$ and $x_2$ is upper triangular. Then we also have another such equation linking $y_1,y_2,y_1',y_2',y_1'',y_2''$ obtained by differentiating  \eqref{eq:deriv-of-j}. When $c\neq 0$, both of these equations depend on $x_1$, so together they yield a single algebraic relation between $y_1,y_2,y_1',y_2',y_1'',y_2''$.
\end{itemize}

\end{remark}

\begin{definition}
For a variety $V \subseteq \mathbb{C}^{3n}$ we let the $J$-\emph{atypical set} of $V$, denoted $\Atyp_J(V)$, be the union of all atypical components of intersections $V \cap T$ in $\mathbb{C}^{3n}$ where $T\subseteq \mathbb{C}^{3n}$ is a $J$-special variety.
\end{definition}

\begin{conjecture}[Modular Zilber-Pink with Derivatives, MZPD {\cite{Pila-MZPD}}]\label{conj:MZPD-Pila}
For every algebraic variety $V \subseteq \mathbb{C}^{3n}$ there is a finite collection $\Sigma$ of proper $\glpq$-special subvarieties of $\mathbb{H}^n$ such that $$\Atyp_J(V) \cap J(\mathbb{H}^n) \subseteq \bigcup_{\substack{U \in \Sigma\\ \bar{\gamma}\in \SL_2(\mathbb{Z})^n}} \langle \langle \bar{\gamma} U \rangle \rangle.$$ 
\end{conjecture}

\begin{remark}
    \begin{itemize}[leftmargin = *]
        \item[]

        \item One could propose a stronger conjecture stating that $\Atyp_J(V)$ is covered by $J$-special varieties corresponding to $\SL_2(\Z)$-translates of finitely many $\glpq$-special varieties $U \in \Sigma$. However, an intersection of $V$ with a $J$-special variety may have a component which does not intersect the image of $J$, or that intersection is small. So while this stronger statement seems sensible (meaning there does not seem to be a trivial counterexample), it is less natural and less about the function $J$ than Conjecture~\ref{conj:MZPD-Pila}. Zilber's original motivation for CIT came from the idea of deducing a uniform version of Schanuel from itself. Similarly, in \cite{Pila-MZPD} Pila proposes MZPD as the difference between MSCD and its uniform version. Since MSCD is about the function $J$, Pila only needed to deal with the part of $\Atyp_J(V)$ that consists of points from the image of $J$. Furthermore, Conjecture~\ref{conj:MZPD-Pila} is supported by the theorems presented in \S\ref{subsec:J-Der-Functional-variants}, while we do not have any evidence towards the said stronger statement, so we do not propose such a conjecture.

        \item Given a $J$-special variety $S \seq \C^{3n}$ with an atypical intersection $V\cap S$, the intersection $\pi_{\ybar}(V) \cap \pi_{\ybar}(S)$ may or may not be atypical. The novelty of MZPD is when this intersection is typical as the atypical ones are accounted for MZP.

        \item In MZPD we may need infinitely many $J$-special varieties to cover the set $\Atyp_J(V) \cap J(\mathbb{H}^n)$ but the conjecture states that it is ``generated'' by finitely many $\GL_2^+(\Q)$-special varieties. See the example below.
    \end{itemize}
\end{remark}

\begin{example}
Consider the variety $V \seq \C^{9}$ defined by a single equation $\Phi_2(y_1,y_2)+\Phi_3(y_2,y_3)=0$. Let $T \seq \C^3$ be a $j$-special variety defined by $\Phi_2(y_1,y_2)=\Phi_3(y_2,y_3)=0$, and let $U \seq \h^3$ be $\glpq$-special such that $j(U)=T$. Then for every $\bar{\gamma} \in (\SL_2(\Z))^3$ we have $\dangle{\bar{\gamma} U}\seq V$, and these are maximal $J$-special (hence atypical) in $V$. Thus, the single $j$-special variety $T$ ``generates'' an infinite set of maximal $J$-atypical subvariaties of $V$.
\end{example}

MZPD has an analytic component -- the intersection of $\Atyp_J(V)$ with the image of $J$. We now propose an ``algebraic'' MZPD conjecture which we believe will be more amenable to (differential) algebraic and geometric techniques (below we provide evidence in support of this). The idea is to replace the set of points from the image of $J$ in an atypical subvariety of $V$ by its Zariski closure. Then we need to understand which algebraic varieties can contain a Zariski dense set of such points, hence this is a variant of the Existential Closedness problem for $J$. So we define an appropriate notion of froadness which serves that purpose.

\begin{definition}
An irreducible variety $W \seq \C^{3n}$ is called $\im(J)$-\emph{froad} (resp. $\im(J)^*$-\textit{froad})\footnote{Here $\im$ stands for the image of a function.} if it is the projection of a $\Gamma_J$-froad (resp. $\Gamma_J^*$-froad) variety $V \seq \h^n \times \C^{3n}$ to the coordinates $(\ybar, \ybar', \ybar'')$.
\end{definition}

 The following statement gives an explicit definition of these notions. Its proof is fairly straightforward from the definitions and is left to the reader.

\begin{proposition}
   Let $W \seq \C^{3n}$ be an irreducible variety, and let $T \seq \C^n$ be the $j$-special closure of $\pi_{\ybar}(W)$. Then $W$ is $\im(J)^*$-froad if and only if there is a $\GL_2^+(\Q)$-special variety $U\seq \h^n$ such that

   \begin{itemize}[leftmargin = *]

    \item  $j(U)=T$,

    \item $W\seq \dangle{U}$,

    \item for any $\bar{i}$ we have
   $\dim \pi_{\bar{i}} (W) \geq \dim \dangle{\pr_{\bar{i}} (U)} - \dim \pr_{\bar{i}} (U). $
   \end{itemize}
 Furthermore, $W$ is $\im(J)$-froad if and only if  $U=\h^n$, $T=\C^n$, and for any $\bar{i}$ of length $k$ we have    $\dim \pi_{\bar{i}} (W) \geq 2k$.

\end{proposition}

\begin{definition}
For a variety $V \subseteq \mathbb{C}^{3n}$ we let the \emph{froadly $J$-atypical set} of $V$, denoted $\FAtyp_J(V)$, be the union of all $\im(J)^*$-froad and atypical components of intersections $V \cap T$ in $\mathbb{C}^{3n}$ where $T\subseteq \mathbb{C}^{3n}$ is a $J$-special variety.
\end{definition}

\begin{conjecture}[Modular Zilber-Pink with Derivatives for Froad varieties, MZPDF]\label{conj:MZPD-broad}
For every algebraic variety $V \subseteq \mathbb{C}^{3n}$ there is a finite collection $\Sigma$ of proper $\glpq$-special subvarieties of $\mathbb{H}^n$ such that 
\[\FAtyp_J(V) \subseteq \bigcup_{\substack{U \in \Sigma\\ \bar{\gamma}\in \SL_2(\mathbb{Z})^n}} \langle \langle \bar{\gamma} U \rangle \rangle.\]
\end{conjecture}

 Now we aim to understand the relation between Conjectures \ref{conj:MZPD-Pila} and \ref{conj:MZPD-broad}. We can show they are equivalent assuming some weakened versions of MSCD and MECD referring only to the image of $J$. We call these conjectures MSCDI and MECDI where ``I'' stands for ``Image''.

\begin{conjecture}[MSCDI]
    Let $z_1,\ldots,z_n \in \h$ be arbitrary and let $U \seq \h^n$ be the $\GL_2^+(\Q)$-special closure of $(z_1,\ldots,z_n)$. Then $$\td_{\Q}\Q(J(z_1),\ldots,J(z_n))\geq \dim \dangle{U} - \dim U.$$
\end{conjecture}

\begin{conjecture}[MECDI]
       Let $V \subseteq \C^{3n}$ be an irreducible $\im(J)^*$-froad variety. Then $V\cap \im(J) \neq \emptyset$.
\end{conjecture}

\begin{proposition}\label{prop:MZPD-equiv}
    \begin{itemize}[leftmargin = 7mm]
        \item[]

        \item[\rm{(i)}] Assume MECDI. Conjecture \ref{conj:MZPD-Pila} (MZPD) implies Conjecture \ref{conj:MZPD-broad} (MZPDF).

        \item[\rm{(ii)}] Assume MSCDI. Conjecture \ref{conj:MZPD-broad} (MZPDF) implies Conjecture \ref{conj:MZPD-Pila} (MZPD).
    \end{itemize}
\end{proposition}

\begin{proof}
    (i) Let $W$ be an $\im(J)^*$-froad atypical subvariety of $V \seq \C^{3n}$. Then by MECDI and  the Rabinowitsch trick (see \cite[Proposition 4.34]{Aslanyan-adequate-predim}), $W\cap \im(J)$ is Zariski dense in $W$. By MZPD (Conjecture \ref{conj:MZPD-Pila}), $\pi_{\ybar} (W \cap \im(J))$ is contained in a union of finitely many $j$-special varieties depending only on $V$. Hence, $\pi_{\ybar} (W) =\overline{\pi_{\ybar} (W \cap \im(J))}^{\Zcl}$ is also contained in that union. Since $W$ is irreducible, $\pi_{\ybar} (W)$ is contained in one such $j$-special variety $T$, and since $W$ is $\im(J)^*$-froad, it is contained in a $J$-special variety associated to $T$.

    (ii) (cf. \cite[Proposition 9.10]{Aslanyan-weakMZPD}) Now assume MSCDI and MZPDF hold. We also assume first that $V$ is defined over $\Qalg$. Let $\bar{w} := (j(\zbar), j'(\zbar), j''(\zbar) ) \in \Atyp_J(V)$ belong to an atypical component of an intersection $V\cap T$ where $T$ is $J$-special. If $T' \seq T$ is the $J$-special closure of $\bar{w}$ (that is, $T' = \dangle{U}$ where $U$ is the $\glpq$-special closure of $\zbar$), then by \cite[Lemma 9.9]{Aslanyan-weakMZPD}, $\bar{w}$ belongs to an atypical component $W$ of the intersection $V \cap T'$. MSCDI implies that $W$ is $\im(J)^*$-froad. Hence, by MZPDF $W$ is contained in a $J$-special variety $S$ associated to one of the finitely many $j$-special varieties depending only on $V$. Then $\bar{w}$ also belongs to $S$.

    When $V$ is defined over arbitrary parameters, rather than $\Qalg$, the same proof goes through provided that we can extend MSCDI and get a lower bound on the transcendence degree of a $J$-point over finitely generated fields. This has been done in \cite[\S 5]{Aslanyan-Eterovic-Kirby-closure-op} for MSCD (see also \cite[\S 4.2]{eterovic-generic}), and MSCDI can be treated similarly.    
\end{proof}

MSCDI, like full MSCD, seems to be out of reach. Hence the second part of the above proposition is not very helpful. On the other hand, MECDI is within reach, albeit still open. 
Therefore, the first part of the proposition is more meaningful and tells us that MZPDF (Conjecture \ref{conj:MZPD-broad}) is probably more tractable than MZPD (Conjecture \ref{conj:MZPD-Pila}). It is unlikely that the second implication in Proposition \ref{prop:MZPD-equiv} can be proven without assuming MSCDI.

\subsection{Functional/differential variants}\label{subsec:J-Der-Functional-variants}

The functional variants of all the above conjectures were established in the last decade. We present them below.

\begin{definition}
Let $(F; +, \cdot, D_1,\ldots,D_m)$ be a differential field with constant field $C = \bigcap_{k=1}^m \ker D_k$. Let also $\Psi$ be the rational function appearing in the differential equation of the $j$-function (see \S \ref{subsec:j-diff-versions}).
\begin{itemize}[leftmargin = *]
    \item We define a $4$-ary relation $\D \Gamma_J(x, y, y', y'')$  by
$$ \exists y''' \left[ \Psi(y,y',y'',y''')=0 \wedge \bigwedge_{k=1}^m D_ky=y'D_kx \wedge D_ky'=y''D_kx \wedge D_ky''=y'''D_kx\right].$$

\item The relation $\D \Gamma_J^{\times}(x,y,y',y'')$ is defined by the formula
$$\D \Gamma_J(x,y,y',y'')\wedge x \notin C \wedge y \notin C \wedge y' \notin C \wedge y'' \notin C.$$

\item The relations $\D \im(J)$ and $\D \im(J)^{\times}$ are defined respectively as $\exists x \D \Gamma_J(x,y,y',y'')$ and $\exists x \D \Gamma_J^{\times}(x,y,y',y'')$.

\item  By abuse of notation, we will use the above expressions ($\D \Gamma_J,\D \Gamma_J^{\times}$,etc.) to denote their Cartesian powers too.
\end{itemize}
\end{definition}

 If $F$ is a field of meromorphic functions of variables $t_1,\ldots,t_m$ over some complex domain with derivations $\frac{d}{d t_k}$, then $\D \Gamma_J^{\times}(F)$ is interpreted as the set of all tuples $(x,y,y',y'')\in F^4$ where $x=x(t_1,\ldots,t_m)$ is some meromorphic function, $y=j(gx)$ for some $g\in \GL_2(\C)$, and $y' = \frac{dj(gx)}{dx}$, $y'' = \frac{d^2j(gx)}{dx^2}$.

The Ax-Schanuel theorem for $J$ is due to Pila and Tsimerman. Again, their proof is based on o-minimality, and Bl\'azquez-Sanz, Casale, Freitag, Nagloo give a differential-algebraic/model-theoretic proof in \cite{sanz-cas-frei-nag-Ax-I}.

\begin{theorem}[Ax-Schanuel for $J$ {\cite[Thm. 1.3]{Pila-Tsim-Ax-j}}]
Let $(F;+,\cdot,D_1,\ldots,D_m)$ be a differential field with commuting derivations and with field of constants $C$. Let also $(z_i, j_i, j_i', j_i'') \in \D \Gamma_J^{\times}(F),~ i=1,\ldots,n.$ If the $j_i$'s are pairwise modularly independent then $\td_CC\left(\bar{z},\bar{j}, \bar{j}', \bar{j}''\right) \geq 3n+\rk (D_kz_i)_{i,k}.$
\end{theorem}

As in the previous section, Ax-Schanuel can be used to prove a differential analogue of MECD.

\begin{theorem}[Differential MECD {\cite[Thm. 1.2]{Aslanyan-Eterovic-Kirby-Diff-EC-j}}]\label{thm-diff-ec-J}
Let $F$ be a differential field, and $V\subseteq F^{4n}$ be a $\Gamma_J$-broad variety. Then there is a differential field extension $K\supseteq F$ such that $V(K)\cap \D \Gamma_J(K) \neq \emptyset$. In particular, if $F$ is differentially closed then $V(F)\cap \D \Gamma_J(F)\neq \emptyset$.
\end{theorem}

\begin{remark}
  In this theorem, when $V$ is defined over the constants $C$ and is \textit{strongly $\Gamma_J$-broad} (i.e. strict inequalities hold in Definition \ref{defin: J broad-free} (first bullet point)), we have $V(K)\cap \D\Gamma_j^{\times}(K) \neq \emptyset$; see \cite[Theorem 1.3]{Aslanyan-Eterovic-Kirby-Diff-EC-j}.
\end{remark}

At the end we state several analogues of MZPD and MZPDF.

\begin{definition}
For a $J$-special variety $T\subseteq \mathbb{C}^{3n}$ and an algebraic variety $V \subseteq \mathbb{C}^{3n}$ an atypical component $W$ of an intersection $V \cap T$ in $\mathbb{C}^{3n}$ is a \emph{strongly $J$-atypical} subvariety of $V$ if for every irreducible analytic component $W_0$ of $W \cap J(\mathbb{H}^n)$, no coordinate is constant on $\pi_{\ybar}(W_0)$. The \emph{strongly $J$-atypical set} of $V$, denoted $\SAtyp_J(V)$, is the union of all strongly $J$-atypical subvarieties of $V$.
\end{definition}

The following is a weak version of MZPD, the proof of which is based on complex geometric tools. It generalises Functional MZP (Theorem \ref{thm:func-MZP}), hence it gives a third proof of the latter.

\begin{theorem}[Weak MZPD {\cite[Thm. 7.9]{Aslanyan-weakMZPD}}]\label{weak-MZPD1-intro}
For every algebraic variety $V \subseteq \mathbb{C}^{3n}$ there is a finite collection $\Sigma$ of proper $\glpq$-special subvarieties of $\mathbb{H}^n$ such that 
$$\SAtyp_J(V) \cap J(\mathbb{H}^n) \subseteq \bigcup_{\substack{U \in \Sigma\\ \bar{\gamma}\in \SL_2(\mathbb{Z})^n}} \langle \langle \bar{\gamma} U \rangle \rangle.$$ 
\end{theorem}

In order to present differential analogues of MZPD(F), we need to introduce several definitions and pieces of notation.

\begin{definition}[{\cite[\S 6]{Aslanyan-weakMZPD}}]
Let $C$ be an algebraically closed field. Define $D$ as the zero derivation on $C$ and extend $(C;+,\cdot, D)$ to a differentially closed field $(K;+,\cdot, D)$. 
\begin{itemize}[leftmargin = *] 
    \item Let $T \subseteq C^n$ be a $j$-special variety and $U\subseteq C^n$ be a $\GL_2(C)$-special variety associated to $T$, that is, $U$ is defined by $\GL_2(C)$-equations and for any $i, k$ the pair of coordinates $x_i, x_k$ are related on $U$ if and only if $y_i, y_k$ are modularly related on $T$. Denote by $\langle \langle U, T \rangle \rangle$ the Zariski closure over $C$ of the projection of the set $$\D \Gamma_J^{\times}(K) \cap (U(K)\times T(K) \times K^2)$$ to the coordinates $(\ybar, \ybar', \ybar'')$.
    
    \item A D$_J$-\emph{special} variety is a variety $S:=\langle \langle U, T \rangle \rangle$ for some $T$ and $U$ as above. 
    
    \item $S \sim T$ means that $S:=\langle \langle U, T \rangle \rangle$ for some $U$ associated to $T$. For a set $\Sigma$ of $j$-special varieties $S \sim \Sigma$ means that $S\sim T$ for some $T \in \Sigma$.
\end{itemize}
\end{definition}

\begin{definition}

 Let $V\subseteq C^{3n}$ be a variety. The D$_J$-\emph{atypical set} of $V$, denoted $\Atyp_{\D_J}(V)$, is the union of all D$_J$-atypical subvarieties of $V$, that is, atypical components of intersections $V \cap T$ where $T\subseteq C^{3n}$ is D$_J$-special. 
 The set $\SFAtyp_{\D_J}(V)$ denotes the union of all $\D_J$-atypical subvarieties of $V$ which are strongly $\im(J)^*$-froad.\footnote{A variety $W\seq \C^{3n}$ is \textit{strongly} $\im(J)^*$-\textit{froad} if there is a $\GL_2^+(\Q)$-special variety $U\seq \h^n$ such that $j(U)=T$, $W\seq \dangle{U}$, and for any $\bar{i}$ we have
   $\dim \pi_{\bar{i}} (W) > \dim \dangle{\pr_{\bar{i}} (U)} - \dim \pr_{\bar{i}} (U).$}
\end{definition}

\begin{theorem}[Functional MZPD; FMZPD {\cite[Thm. 8.2]{Aslanyan-weakMZPD}}]\label{DMZPD-intro}
Let $(K; +,\cdot, D)$ be a differential field with an algebraically closed field of constants $C$. Given an algebraic variety $V\subseteq C^{3n}$, there is a finite collection $\Sigma$ of proper $j$-special subvarieties of $C^n$ such that
\[\Atyp_{\D_J}(V)(K) \cap \D \im_J^{\times}(K) \subseteq \bigcup_{\substack{S\sim \Sigma}} S.\]
\end{theorem}

\begin{theorem}[Functional MZPDF; FMZPDF {\cite[Thm. 9.8]{Aslanyan-weakMZPD}}]\label{FMZPD-intro}
Let $C$ be an algebraically closed field of characteristic zero. Given an algebraic variety $V\subseteq C^{3n}$, there is a finite collection $\Sigma$ of proper $j$-special subvarieties of $C^n$ such that 
\[\SFAtyp_{\D_J}(V)(C) \subseteq \bigcup_{\substack{S\sim \Sigma}} S.\]
\end{theorem}

These theorems are analogues of MZPD and MZPDF respectively and so they support those conjectures. In \cite{Aslanyan-weakMZPD} we give a complex geometric proof of FMZPD (the transition from complex geometry to differential algebra is via Seidenberg's embedding theorem) and a differential-algebraic proof of FMZPDF. The core of both proofs is the Ax-Schanuel theorem for $J$. The proof of FMZPDF also uses the differential version of MECDI which is a special case of Theorem~\ref{thm-diff-ec-J}. As above, FMZPD and FMZPDF can be deduced from one another using Differential MECDI, so that gives two proofs for each of the above theorems, one differential algebraic and one complex geometric.

For further results on MECD and MZPD see \cite{Aslanyan-Kirby-blurred-j,Aslanyan-Eterovic-Kirby-closure-op,eterovic-generic} and \cite{Aslanyan-weakMZPD} respectively. Spence \cite{Spence} has proven some results towards the Modular Andr\'e-Oort with Derivatives conjecture which is a special case of MZPD.

\begin{remark}
     Section \ref{sec: incorporating derivatives} turned out to be somewhat technical with some hard-to-remember concepts and notation. Unfortunately, that seems necessary for precision and rigour. A reader, who is not familiar with the general topics discussed here, may be lost in the various versions of the conjectures and theorems. Therefore, we would like to reiterate the main high-level idea of this section: incorporating the derivatives into the modular versions of the Schanuel, Existential Closedness, and Zilber-Pink conjectures gives a deeper insight into these problems and reveals some hidden (possibly surprising) connections between them. Exploring these conjectures in this more general setting would allow us to better understand the full model-theoretic picture.
\end{remark}

\addtocontents{toc}{\protect\setcounter{tocdepth}{1}}
\subsection*{Acknowledgement.} I would like to thank the anonymous referees for their useful comments.
\addtocontents{toc}{\protect\setcounter{tocdepth}{2}}

\bibliographystyle {alpha}
\bibliography {ref}

\end{document}